\documentclass[11pt,twoside]{article}
\usepackage{amsfonts,amsmath,epsfig,amssymb,color,amsthm}

\usepackage{cite}

\usepackage{xcolor}

\numberwithin{equation}{section}
\newfont{\ctv}{msam10}
\newcommand{\bbox}{\mbox{\ctv \symbol{3}}}
\def\QED{{$\hfill\bbox$}}
\newenvironment{pf}[1]{\par\vskip1mm{\noindent\it #1.}\ }{\QED\par\vskip2mm}

\newtheorem{proposition}{Proposition}[section]
\newtheorem{lemma}{Lemma}[section]
\newtheorem{definition}{Definition}[section]

\newtheorem*{theoremH2.2}{Hypothesis 2.2}

\newcommand{\Barint}{-\kern-.15in\int}
\newcommand{\barint}{-\kern-.12in\int}
\DeclareMathOperator{\domain}{dom}

\def\bpf{\begin{pf}}
\def\epf{\end{pf}}

\begin{document}

\thispagestyle{empty}
 \setcounter{page}{1}

\title{Periodic solutions of a phase-field model with hysteresis \footnote{Supported by National Natural Science Foundation of China (12071165 and 62076104), Natural Science Foundation of Fujian Province (2020J01072), Program for Innovative Research Team in Science and Technology in Fujian Province University, Quanzhou High-Level Talents Support Plan (2017ZT012), Scientific Research Fund of Huaqiao University (605-50Y19017 and 605-50Y14040). The research of the second author was also supported in part by RFBR grant no. 18-01-00026}
}
 \author{Chen Bin\footnote{Fujian Province University Key Laboratory of Computational Science, School of Mathematical Sciences, Huaqiao University, Quanzhou 362021, China, E-mail: {\tt
chenbinmath@163.com}. }
\ and Sergey A. Timoshin\footnote{Fujian Province University Key Laboratory of Computational Science, School of Mathematical Sciences, Huaqiao University, Quanzhou 362021, China -- and --
Matrosov Institute for System Dynamics and Control
 Theory,
 Russian Academy of Sciences, Lermontov str. 134,
664033  Irkutsk,
 Russia, E-mail: {\tt
sergey.timoshin@gmail.com} (Corresponding
author).} }

\date{\today}

\maketitle \hyphenation{Ca-ra-the-odo-ry}

\begin{abstract}
\noindent {

In the present paper we consider a partial differential system describing a phase-field model with temperature dependent constraint for the order parameter. The system consists of
an energy balance equation with a fairly general
nonlinear heat source term and a phase dynamics
equation which takes into account the hysteretic character of the
process. The existence of a periodic solution for this system is proved under a minimal set of assumptions on the curves defining the corresponding hysteresis region.
}
\end{abstract}





\section{Introduction} \label{number of section:introduction}

In the space-time cylinder $Q := [0,T]\times\Omega$, where
$\Omega \subset \mathbb{R}^N (N\geq 1)$ is a bounded domain with
smooth boundary  $\partial\Omega$ and $T>0$ is a fixed final time,
consider the system

\begin{equation}
 u_t - \Delta u =h(u,v) \quad \text{in} \; Q
,\label{1.1}
\end{equation}
\begin{equation} v_t -\kappa \Delta v+ \partial I(u;v) \ni g(u,v)  \quad \text{in}
\; Q , \label{1.2}\end{equation}
\begin{equation}u =v= 0 \quad \text{on}
\; (0,T)\times\partial\Omega , \label{1.3}\end{equation}
\begin{equation}u(x,0) = u(x,T)
, \quad v(x,0) = v(x,T) \quad \text{on} \; \Omega .
\label{1.4}\end{equation}

Here,  $I(u;\cdot)$ is the
indicator function of the interval $[f_*(u), f^*(u)]$, $\partial
I(u;\cdot)$ is its subdifferential in the sense of convex analysis,  $h, g, f_*, f^*$ are given
functions with the properties specified in the next section, $\kappa>0$ is a given constant.

\hspace{0.6cm}  For convenience, denote system
$(\ref{1.1})$--$(\ref{1.4})$ by $(P)$. System $(P)$ can be regarded as a dynamical model of a phase transition process between two distinct phases (such as solid-liquid) placed in the container $\Omega$.
The state variables $u=u(t,x)$
and $v=v(t,x)$ are then interpreted as the relative temperature and the order parameter (phase fraction of an individual phase), respectively. Eq. $(\ref{1.2})$ with $g\equiv 0$ and $\kappa=0$
models a continuous hysteresis operator of generalized play type generated by the
curves $v=f_*(u)$ and  $v=f^*(u)$, see \cite{BS, Krejci, Visintin}
for details. The introduction of the latter operator to the model
accounts for hysteretic relationship between $u$ and $v$, playing in
this case the roles of the input and output functions, respectively.

\hspace{0.6cm} Recent years have seen a considerable amount of works on partial differential equations with hysteresis. In particular, the questions on existence, uniqueness and large time behaviour of solutions to Cauchy problems for systems with state-dependent constraints,  relevant to $(\ref{1.1})$--$(\ref{1.3})$ were addressed in a number of papers (see, e.g., \cite{KKV, KeV, KeKM, CKK, KMO, Kubo, KTT, AT}, and references therein). Periodic problems for systems with hysteresis describing phase transitions have also received a keen attention. Among related contributions, we mention the works \cite{Krejci1, Krejci2, Krejci3} studying periodic processes in elastoplastic bodies, \cite{GHM} dealing with  the Stefan problem in a one-dimensional spatial domain, that involves a relay hysteresis operator, and \cite{comwa} which considers the system $(\ref{1.1})$--$(\ref{1.4})$ with $h(u,v)=v$, $g\equiv 0$, and $\kappa=0$. We note, however, that the requirements on the functions $f_*$ and $f^*$ in \cite{comwa}: $f_*, f^*\in C^2(\mathbb{R})\cap L^\infty(\mathbb{R})$ are non-decreasing Lipschitz continuous, $f_*(u)=f^*(u)$ for $u\in (-\infty, a]\cup[b, +\infty)$ with $a<0<b$, $f_*$ is convex on $(-\infty, b)$ and $f^*$ is concave on $(a, +\infty)$, which are indispensable for the proof in \cite{comwa}, appear to be too demanding. 
In our paper, we dispense with these assumptions on the functions $f_*$ and $f^*$ retaining only the Lipschitz continuity. In this respect, as a byproduct of our analysis, we also improve the results on the existence of a solution to the Cauchy problem from \cite{CKK} by removing the assumptions of smoothness, monotonicity and boundedness on the functions $f_*$ and $f^*$. Moreover, the convergences of approximate solutions which the authors obtain in \cite[The proof of Theorem 2.1]{comwa} would not be sufficient to treat general nonlinear right-hand sides as in $(\ref{1.1})$, $(\ref{1.2})$. Another advance of our paper is that we allow a diffusion effect for $v$ assuming the coefficient of the interfacial energy $\kappa$ to be nonzero.

\hspace{0.6cm} Here, we would like to mention that hysteresis curves encountered in the practice of physical measurements, stress-strain hysteresis loops in shape memory alloy wires, load-displacement hysteresis curves in composite structures, magnetic hysteresis curves of nano-minerals may genuinely occur to lack the smoothness.  Also, we add that for the completed relay operator and the truncated play operators employed to approximate the former (see \cite{Visintin})  the curves describing the corresponding hysteresis regions are piecewise linear but nonsmooth.

\hspace{0.6cm} The purpose of the present paper is to prove the existence of a solution to system $(P)$ with sufficiently general $h$, $g$, $f_*$, and $f^*$. In our approach to establish the existence for
problem $(P)$, first, we construct  a family of suitable approximate problems based on the Yosida
regularization $\partial I^\lambda(u;\cdot)$, $\lambda>0$, of
the subdifferential $\partial I(u;\cdot)$. We then further regularize the nonsmooth functions $f_*$ and $ f^*$ describing $\partial I^\lambda(u;\cdot)$ by  sequences of mollifiers depending on a regularizing parameter $\varepsilon>0$. Next, we rewrite thus obtained approximate problem $(P)_{\lambda,\varepsilon}$ as a single abstract differential equation with a general nonlinearity subject to periodic condition. Applying a standard for such equations technique invoking a fixed point argument we find a solution of the abstract equation which provides us with a solution to the approximate problem $(P)_{\lambda,\varepsilon}$.
After that, we establish a priori estimates independent of
$\varepsilon$ for solutions of the approximate problems and performing a limiting procedure as $\varepsilon\to 0$ we obtain an intermediate approximate problem $(P)_{\lambda}$ depending now on the parameter $\lambda$ only. Deriving uniform estimates with respect to $\lambda$ for the latter system, we finally
prove the existence of a solution to problem $(P)$ through the
passage-to-the-limit procedure when $\lambda\to 0$. We note that in order to get
suitable compactness properties and, thus, to
legitimate this passage-to-the-limit we exploit essentially the
properties derived from the specific structure of the
approximate equations for $(\ref{1.2})$. This also allows us,
inter alia, to treat general nonlinearities in Eqs.
$(\ref{1.1})$, $(\ref{1.2})$.

\section{Preliminaries and hypotheses on the data} \label{section2}

In this section, we recall some notions which we use in the paper and posit assumptions on the
data describing  Problem $(P)$.

\hspace{0.6cm} Throughout the paper, we denote by $H$ the Hilbert
space $L^2(\Omega)$ with the standard inner product
$\langle\cdot,\cdot \rangle_H$, and by $V$ the Sobolev space $H^1(\Omega)$.

\hspace{0.6cm} Given a Hilbert
space $X$ with the inner product
$\langle\cdot,\cdot \rangle_X$ and  a convex, lower semicontinuous function $\varphi:
X \to {\mathbb R}\cup \{+\infty\}$ which is not identically
$+\infty$, the subdifferential   $\partial \varphi (x)$ of $\varphi$  at a
point $x\in X$ is, in general, a set  defined by the rule
$$
\partial\varphi (x)=\{h\in X; \; \langle h,y-x\rangle_X\leq
\varphi(y)-\varphi(x) \;\; \forall y\in X\}.
$$

\hspace{0.6cm} Let $f_*, f^*$ be two Lipschitz continuous functions
defined on $\mathbb{R}$. Then, the subdifferential of the
indicator function $I(u;\cdot)$, $u\in \mathbb{R}$,
\begin{equation*} I(u;v):=  \left\{
\begin{array}{cl}
0 &\mbox{if} \;\; \; f_{*}(u)\leq v\leq  f^{*}(u) ,  \\
 +\infty &\mbox{otherwise} ,
\end{array}
\right.
\end{equation*}  of the
interval $[f_{*}(u), f^{*}(u)]$ has the form:

\begin{equation}\partial I(u;v)= \left\{ \begin{array}{ccl}
\emptyset &\mbox{if}& v \notin [f_{\ast}(u),f^{\ast}(u)] ,\\
{[}0,+\infty) &\mbox{if}& v = f^{\ast}(u)> f_{\ast}(u),\\
\{0\}       &\mbox{if}& f_{\ast}(u)< v < f^{\ast}(u) ,\\
(-\infty, 0{]} &\mbox{if}& v = f_{\ast}(u)< f^{\ast}(u),\\
(-\infty, +\infty) &\mbox{if}& v = f_{\ast}(u)= f^{\ast}(u).
\end{array}
\right. \label{2.1}
\end{equation}

For $\lambda>0$, the Yosida regularization of $\partial I(u;v)$ is
the function
\begin{equation}\partial I^\lambda(u;v)= \frac{1}{\lambda}[v-f^*(u)]^+-\frac{1}{\lambda}[f_*(u)-v]^+ , \quad u,v\in \mathbb{R}.\label{2.2}
\end{equation}

\hspace{0.6cm} Let $f:\mathbb{R}\to \mathbb{R}$ be a Lipschitz continuous function. For $\varepsilon>0$, denote by  $f_{\varepsilon}(u)$, $u\in \mathbb{R}$,  the following
regularization of the function $f_{\varepsilon}(u)$:

\begin{align}
f_\varepsilon(u):=\int_{\mathbb{R}}f(s)\rho_\varepsilon(u-s)\, ds =\int_{\mathbb{R}}f(u-\varepsilon s)\rho(s)\, ds , \label{2.3}
\end{align}

where $\rho\in C^\infty(\mathbb{R})$ is such that $\rho\geq 0$,
$\rho(s)=0$ when $|s|\geq 1$, $\rho(s)=\rho(-s)$,
$\int_{\mathbb{R}}\rho(s) \,ds=1$,
$\rho_\varepsilon(s):=\varepsilon^{-1}
\rho\left(\frac{s}{\varepsilon}\right)$.

\hspace{0.6cm} The lemma below follows directly from the
definition of $f_\varepsilon(u)$ and the properties of
$\rho_\varepsilon(s)$.

\begin{lemma} \label{lemma2.1} The function $f_\varepsilon(u)$
possesses the following properties:
\begin{itemize}
\item[{$(1)$}] $f_\varepsilon(u)\in C^\infty(\mathbb{R});$

\item[{$(2)$}] $f_\varepsilon(u)$ is Lipschitz continuous with the
same Lipschitz constant as $f(u);$

\item[{$(3)$}] $f_\varepsilon(u)\to f(u)$ as $\varepsilon\to 0$
uniformly on $\mathbb{R};$

\item[{$(4)$}] $|f_\varepsilon|_\infty \leq |f|_\infty.$
\end{itemize}
\end{lemma}

\hspace{0.6cm} Problem $(\ref{1.1})$--$(\ref{1.4})$ is considered
under the following hypotheses:

\begin{itemize}

\item[{\textbf{(H1)}}] the functions $h, g:\mathbb{R}^2\to \mathbb{R}$ are locally Lipschitz continuous, the function $h$ is bounded, and the Lipschitz constant in the variable $v$ of $g$: $L_*<\kappa/C_P$, where $C_P$ is the best constant of the Poincar\'{e} inequality on $\Omega$;

\item[{\textbf{(H2)}}] the functions $f_*, f^*:\mathbb{R}\to \mathbb{R} $ are locally Lipschitz continuous, $f_*(u)\leq  f^*(u)$ for all $u\in \mathbb{R}$, and there exist constants $a, b$, $a<b$, such that $f_*(u)= f^*(u)$ for $u\in \mathbb{R}\setminus (a,b)$.

\end{itemize}

\hspace{0.6cm} Next, we define a notion of solution to our Problem
$(P)$.

\begin{definition} \label{def2.1} A pair $\{u,v\}$ is called a solution of
 system $(\ref{1.1})$--$(\ref{1.4})$  if
\begin{itemize}
\item[$(i)$] $ u, v\in W^{1,2}(0,T;H)\cap
L^\infty(0,T;V)\cap L^2(0,T;H^2 (\Omega ))$;

\item[$(ii)$] $u'-\Delta u = h(u,v)$ \; in $H$ a.e. on $[0,T]$;

\item[$(iii)$] $ v'-\kappa \Delta v+ \partial I(u;v)\ni g(u,v)$
\; in $H$ a.e. on $[0,T]$;

\item[$(iv)$] $u=v=0$ on $\partial\Omega$ $($in the sense of traces$)$ a.e. on $[0,T]$;

\item[$(v)$]  $u(0)=u(T)$, $v(0)=v(T)$ \; in $H$,
\end{itemize}

where the prime denotes the derivative with respect to $t$.
\end{definition}

\hspace{0.6cm} We note that in view of the variational characterization of subdifferential, the part $(iii)$ of the above definition implies that
\begin{equation}  f_*(u)\leq v \leq f^*(u)  \;\; \text{a.e. in } \; Q(T)\label{2.4}
\end{equation}
and
\begin{equation}  \langle v'(t)-\kappa \Delta v(t)- g(u(t),v(t)),
    v(t)-z\rangle_H\leq 0\label{2.5}
\end{equation}

for all $z\in H$ with $f_*(u(t))\leq z \leq f^*(u(t))$ a.e. in $\Omega$ for a.e. $t\in [0,T]$.

\section{Approximate problems} \label{section4}

First, we note that if $(u,v)$ is a solution to Problem $(P)$, then $(\ref{1.1})$, $(\ref{1.3})$, $(\ref{1.4})$, and Hypothesis $(H1)$ imply that $|u|_{L^\infty(Q)}\leq M_1$, where $M_1>0$ is a constant depending on $|h|_{L^\infty(\mathbb{R}^2)}$ and $T$ only (see, e.g., \cite[Proposition 10]{Quittner}). By virtue of this estimate, from Hypothesis $(H2)$ and $(\ref{2.4})$ we see that also $|v|_{L^\infty(Q)}\leq M_2$ for a constant $M_2>0$ depending on $M_1$ and the Lipschitz constants of $f_*, f^*$. Therefore, we may now assume (cutting off outside the set where $u$ and $v$ are bounded, if necessary) that the functions $h, g, f_*, f^*$ are all bounded and globally Lipschitz continuous.

\hspace{0.6cm}  In order to prove the existence of a solution to our Problem $(P)$, we approximate the latter by a family of suitable
problems depending on two approximation parameters which we
introduce next.

\hspace{0.6cm} Let $\varepsilon>0$ and ${f}_{*\varepsilon}(u)$ and ${f}^*_\varepsilon(u)$ be the regularizations as in $(\ref{2.3})$ of the functions ${f}_*(u)$ and ${f}^*(u)$, respectively.

\hspace{0.6cm} For $\lambda, \varepsilon>0$, we consider the following
approximate periodic problem denoted by $(P)_{\lambda,\varepsilon}$:
\begin{equation}
 u' -\Delta u =h(u,v) \quad \text{in } H  \; \text{a.e. on } [0,T]
,\label{3.1}
\end{equation}
\begin{equation} v'-\kappa\Delta v+ \partial {I}^{\lambda}_{\varepsilon}(u;v)= g(u,v)   \quad \text{in } H  \; \text{a.e. on } [0,T] , \label{3.2}\end{equation}
\begin{equation}u =  v=0 \quad \text{on } \partial\Omega \;\; \text{a.e. on } [0,T], \label{3.3}\end{equation}
\begin{equation}u(0) =
u(T), \quad v(0) = v(T) \quad \text{in } H ,
\label{3.4}\end{equation}

where $\partial {I}^{\lambda}_{\varepsilon}(u;v)$ is defined as $\partial {I}^{\lambda}(u;v)$ in $(\ref{2.2})$ with ${f}_*$ and ${f}^*$ replaced by ${f}_{*\varepsilon}$ and ${f}^*_\varepsilon$, respectively.

\hspace{0.6cm} A pair of functions $\{u,v\}$ is called a solution to
$(P)_{\lambda,\varepsilon}$ if $u,v\in W^{1,2}(0,T;H)\cap L^\infty(0,T;V)\cap
L^2(0,T;H^2 (\Omega ))$ and  $(\ref{3.1})$--$(\ref{3.4})$ hold.

\hspace{0.6cm} In this section, we prove the existence of solutions for
problems $(P)_{\lambda,\varepsilon}$, $\lambda,\varepsilon>0$. To this aim, first, for convenience, we rewrite    $(P)_{\lambda,\varepsilon}$ as the following equivalent to $(P)_{\lambda,\varepsilon}$ periodic problem for a single abstract nonlinear differential equation in the Hilbert space $\mathbf{H}:=H\times H$:
\begin{equation} z'(t)+ \partial \varphi(z(t))= F(z(t))   \quad \text{in } \mathbf{H}  \quad \text{for a.e.  } t\in [0,T] , \label{3.5}\end{equation}
\begin{equation}z(0) =
z(T) \quad \text{in } \mathbf{H}.
\label{3.6}\end{equation}

Here, $z:=(u,v)\in\mathbf{H}$,
\begin{equation}\varphi(z):= \left\{ \begin{array}{ll}
\displaystyle{ \frac{1}{2}|\nabla u|_H^2+\frac{\kappa}{2}|\nabla v|_H^2 +\frac{1}{2\lambda}|v|_H^2}&\mbox{if } z \in H_0^1(\Omega)\times H_0^1(\Omega) ,\\
+\infty &\mbox{otherwise},
\end{array}
\right. \label{3.7}
\end{equation}
\begin{equation} F(z):= \left( \begin{array}{c}
 h(u,v) \\
g(u,v)+\displaystyle{\frac{1}{\lambda}J_uv}
\end{array}
\right), \label{3.8}
\end{equation}

where $J_uv:=\max\{\min\{v,{f}^*_\varepsilon(u)\},{f}_{*\varepsilon}(u)\}$ is the projection of $v$ onto the set $K(u):=\{v\in H: \; v(x)\in [{f}_{*\varepsilon}(u(x)),{f}^*_\varepsilon(u(x))] \mbox{ for a.e. } x\in \Omega\}$ (cf. $(\ref{2.2})$). Then, $\varphi:\mathbf{H}\to [0,+\infty]$ is a proper, convex, lower semicontinuous function,
\begin{equation} \partial\varphi(z)= \left( \begin{array}{c}
 -\Delta u \\
-\kappa\Delta v+\displaystyle{\frac{1}{\lambda}v}
\end{array}
\right),  \domain \partial\varphi=(H_0^1(\Omega)\times H_0^1(\Omega))\cap (H^2(\Omega)\times H^2(\Omega)), \label{3.9}
\end{equation}

and the notion of a solution to $(\ref{3.5})$, $(\ref{3.6})$ (as well as to $(\ref{3.10})$, $(\ref{3.11})$ below) naturally extends from Definition 2.1.

\hspace{0.6cm}  Given a function $f\in L^2(0,T;\mathbf{H})$ consider the following periodic problem associated with $(\ref{3.5})$, $(\ref{3.6})$:
\begin{equation} z'(t)+ \partial \varphi(z(t))= f(t)   \quad \text{in } \mathbf{H}  \quad \text{for a.e.  } t\in [0,T] , \label{3.10}\end{equation}
\begin{equation}z(0) =
z(T) \quad \text{in } \mathbf{H}.
\label{3.11}\end{equation}

It is a standard matter to show that the function $\varphi$ defined by $(\ref{3.7})$ satisfies the assumptions $(A.1)$--$(A.3)$ of \cite{Yamada}. Moreover, since the subdifferential $\partial\varphi$ (see $(\ref{3.9})$) of this function is obviously strictly monotone operator, \cite[Theorem 1.4]{Yamada} (see, also, \cite[\S 2.3]{Kenmochi}) implies that there exists a unique solution $z=z(f)$ of $(\ref{3.10})$, $(\ref{3.11})$.

\hspace{0.6cm}  Next, define the solution operator $\mathcal{T}:L^2(0,T;\mathbf{H})\to C([0,T];\mathbf{H})$
which with each $f\in L^2(0,T;\mathbf{H})$ associates the unique solution $z=\mathcal{T}(f)$ of
Problem $(\ref{3.10})$, $(\ref{3.11})$.

\begin{proposition} \label{Proposition 3.1}
The operator $\mathcal{T}:L^2(0,T;\mathbf{H})\to C([0,T];\mathbf{H})$
is weak-strong continuous in the sense that if $f_n\to f$ weakly in $L^2(0,T;\mathbf{H})$ for $f, f_n\in L^2(0,T;\mathbf{H})$, $n\geq 1$, then $\mathcal{T}(f_n)\to \mathcal{T}(f)$ strongly in  $C([0,T];\mathbf{H})$.
\end{proposition}

\bpf{Proof} Let $f_n\to f$ weakly in $L^2(0,T;\mathbf{H})$ and $z_n:=\mathcal{T}(f_n)$, $n\geq 1$. Then, we have
\begin{equation} z_n'(t)+ \partial \varphi(z_n(t))= f_n(t)   \quad \text{in } \mathbf{H}  \quad \text{for a.e.  } t\in [0,T] , \label{3.12}\end{equation}
\begin{equation}z_n(0) =
z_n(T) \quad \text{in } \mathbf{H}.
\label{3.13}\end{equation}

Testing Eq. (\ref{3.12}) by $ z'_{n}(t)$ and applying Young's inequality we obtain (cf. $(\ref{3.7})$)
\begin{equation} |z'_{n}|_\mathbf{H}^2+\frac{d}{dt}\langle \partial \varphi(z_n), z_n\rangle_\mathbf{H} \leq C_1 \quad \text{ a.e. on} \;
[0,T],
\label{3.14}\end{equation}

where $C_1:=\sup\limits_{n\geq 1}|f_n|^2_\mathbf{H}$, $|\cdot|_\mathbf{H}$ and $\langle \cdot, \cdot\rangle_\mathbf{H}$ are the norm and inner product in $\mathbf{H}$, respectively. Similarly, testing  Eq. (\ref{3.12}) by $\partial \varphi(z_n)$ we have
\begin{equation} \frac{d}{dt}\langle \partial \varphi(z_n), z_n\rangle_\mathbf{H} +|\partial \varphi(z_n)|_\mathbf{H}^2\leq C_1\quad \text{ a.e. on} \;
[0,T].
\label{3.15}\end{equation}

Now, we test $(\ref{3.12})$ by $z_{n}$ and invoke the Poincar\'{e} inequality  to obtain
\begin{equation*} \frac{d}{dt}|z_{n}|^2_\mathbf{H} + \langle \partial \varphi(z_n), z_n\rangle_\mathbf{H}\leq C_2 \quad  \text{a.e. on} \; [0,T], \label{}\end{equation*}

where $C_2>0$ is a constant depending on $C_1$, $\lambda$ and the constant of the Poincar\'{e} inequality.  Integrating this inequality from $0$ to $T$ and taking account of the periodicity condition $(\ref{3.13})$  we deduce that
\begin{equation} \int_0^T\langle \partial \varphi(z_n(\tau)), z_n(\tau)\rangle_\mathbf{H} d\tau\leq C_2 T.\label{3.16}\end{equation}

Next, to show that the sequence
\begin{equation}  \{\langle \partial \varphi(z_n), z_n\rangle\}_{n\geq 1} \;\;\text{ is bounded in } L^\infty(0,T;\mathbf{H})\label{3.17}\end{equation}

 we note that $z_n$, $n\geq 1$, are $T$-periodic and we can thus consider their periodic extensions onto the whole real axis. Then, we take $t\in[T,2T]$, $\tau\in (0,T)$, and integrate the inequality
\begin{equation*} \frac{d}{dt}\langle \partial \varphi(z_n), z_n\rangle_\mathbf{H}\leq C_1 \quad  \text{a.e. on} \; [0,T], \label{}\end{equation*}

obtained from $(\ref{3.15})$, from $\tau$ to $t$ to get
\begin{equation*} \langle \partial \varphi(z_n(t)), z_n(t)\rangle_\mathbf{H} \leq \langle \partial \varphi(z_n(\tau)), z_n(\tau)\rangle_\mathbf{H} +2C_1 T .\label{}\end{equation*}

Integrating this inequality over $\tau$ from $0$ to $T$ and using $(\ref{3.16})$ we obtain $(\ref{3.17})$.

\hspace{0.6cm} Integrating $(\ref{3.14})$, $(\ref{3.15})$ from $0$ to $T$ and taking account of the periodicity condition $(\ref{3.13})$ we see in view of  $(\ref{3.7})$  that the following uniform  with respect to  $n\geq 1$ estimates hold for the  components $u_n, v_n$ of the solution $z_n$ of $(\ref{3.12})$, $(\ref{3.13})$:
\begin{align} \nonumber
|u_{n}'|_{L^2(0,T;H)}&+|\Delta u_{n}|_{L^2(0,T;H)}
+|\nabla u_{n}|_{L^\infty(0,T;H)}\\+
|v_{n}'|_{L^2(0,T;H)}&+|\Delta v_{n}|_{L^2(0,T;H)}
+|\nabla v_{n}|_{L^\infty(0,T;H)}\leq C_3,\label{3.21}
\end{align}

 for a positive constant $C_3$ independent of $n\geq 1$. On account of these uniform estimates, by weak and
weak-star compactness results, there exists a  subsequence (still indexed by $n$) of the sequence
$z_n$, $n\geq 1$,  and a function $z\in W^{1,2}(0,T; \mathbf{H}) \cap L^{\infty}(0,T;V\times V) \cap
L^2(0,T;H^2(\Omega)\times H^2(\Omega))$ such that
\begin{align}
& \left. \begin{array}{lc} & {z}_{n} \to z \quad
\mbox{weakly in } W^{1,2}(0,T; \mathbf{H}) \cap
L^2(0,T;H^2(\Omega)\times H^2(\Omega))  \vspace{0.1cm}\\
& \*\hspace{3cm} \mbox{and weakly-star in } L^{\infty}(0,T;V\times V).
\end{array}
\right. \label{3.22}
\end{align}

In particular, we also have
\begin{equation} {z}_{n} \to z \quad  \mbox{in } C([0,T];\mathbf{H}).\label{3.23}\end{equation}

The convergences $(\ref{3.22})$, $(\ref{3.23})$ imply that $z$ satisfies $(\ref{3.10})$, $(\ref{3.11})$. Therefore, $z=\mathcal{T}(f)$
and the claim of the proposition follows.
\epf

\hspace{0.6cm}  From Hypotheses $(H1), (H2)$, Lemma 2.1 (2), (4), and $(\ref{3.8})$ we see that the mapping $F:C([0,T];\mathbf{H})\to L^2(0,T;\mathbf{H})$ is continuous and there exists a constant $R>0$ such that
\begin{equation} |F(z)|_\mathbf{H} \leq R \quad  \mbox{for all } z\in \mathbf{H}.\label{3.24}\end{equation}

We now introduce the set
\begin{equation*}
S_R=\{f\in L^2(0,T;\mathbf{H}): \;|f(t)|_\mathbf{H}\leq R
\quad\text{for a.e.} \; t\in [0,T]\}, \label{}
\end{equation*}

and define the superposition $F\circ \mathcal{T}: S_R\to L^2(0,T;\mathbf{H})$ of the solution operator $\mathcal{T}$  and
$F$:
\begin{equation*}
F\circ \mathcal{T}(f)=F(\mathcal{T}(f)). \label{}
\end{equation*}

From Proposition 3.1 and $(\ref{3.24})$ it follows that
$F\circ \mathcal{T}:S_R\to S_R$ is weak-weak continuous.
Since the set $S_R$ is obviously convex and compact in the weak topology of
the space $L^2(0,T;\mathbf{H})$, from the Schauder fixed point theorem we
conclude that there exists a fixed point $f_*\in S_R$ of the
operator $F\circ \mathcal{T}$:
\begin{equation*} f_*=F(\mathcal{T}(f_*)) .\label{}\end{equation*}

Setting $z_*:=\mathcal{T}(f_*)$
we see that $z_*$ is a solution to the periodic problem $(\ref{3.5})$, $(\ref{3.6})$, which provides the desired solution of Problem $(P)_{\lambda, \varepsilon}$, $\lambda, \varepsilon>0$.

\section{Well-posedness of problem $(P)$}

In this section, first we derive uniform a priori estimates independent of the parameter $\varepsilon>0$ for solutions $\{u_{\lambda \varepsilon},v_{\lambda \varepsilon}\}$ of the approximate periodic Problem $(P)_{\lambda, \varepsilon}$, which will allow us to derive the convergence of $\{u_{\lambda \varepsilon},v_{\lambda \varepsilon}\}$ as $\varepsilon\to 0$ to a solution $\{u_{\lambda},v_{\lambda}\}$ of an intermediate approximate problem depending on the parameter $\lambda$ only. Then, we establish uniform bounds independent of the parameter $\lambda$ for solutions $\{u_{\lambda},v_{\lambda}\}$ of the latter system and finally pass to the limit as $\lambda\to 0$ to obtain a solution to our original periodic problem $(P)$.

\subsection{Passage-to-the-limit: $\varepsilon\to 0$}

 Fix $\lambda>0$. Repeating the reasoning in derivation of $(\ref{3.21})$ from $(\ref{3.12})$, $(\ref{3.13})$, using Hypothesis $(H1)$ we obtain from $(\ref{3.1})$, $(\ref{3.3})$  the following uniform  with respect to the parameter $\varepsilon>0$ estimates  for the first component of solutions $(u_{\lambda \varepsilon}, v_{\lambda \varepsilon})$ of the approximate periodic problems $(P)_{\lambda \varepsilon}$:
\begin{align} 
|u_{\lambda \varepsilon}'|_{L^2(0,T;H)}&+|\Delta u_{\lambda \varepsilon}|_{L^2(0,T;H)}
+|\nabla u_{\lambda \varepsilon}|_{L^\infty(0,T;H)}\leq R_1,\label{4.1}
\end{align}

 for a positive constant $R_1$ depending on $|h|_{L^\infty(\mathbb{R}^2)}$ and the Lebesgue measure of $\Omega$, but independent of $\varepsilon$. On account of these uniform estimates, by weak and
weak-star compactness results, there exists a null sequence
$\varepsilon_n$, $n\geq 1$, in $(0,1]$ and a function $u_\lambda$ such that
\begin{align}
& \left. \begin{array}{lc} & {u}_{\lambda \varepsilon_n} \to u_\lambda \quad
\mbox{weakly in } W^{1,2}(0,T; H) \cap
L^2(0,T;H^2(\Omega))  \vspace{0.1cm}\\
& \*\hspace{3cm} \mbox{and weakly-star in } L^{\infty}(0,T; V).
\end{array}
\right. \label{4.2}
\end{align}

In particular, we also have
\begin{equation} {u}_{\lambda \varepsilon_n} \to u_\lambda \quad  \mbox{in } C([0,T];H).\label{4.3}\end{equation}

Invoking the Poincar\'{e} inequality from $(\ref{4.1})$ we obtain
\begin{equation} \int_0^T|u_{\lambda \varepsilon}(\tau)|^2_H d\tau\leq C_P^2 R_1.\label{4.4}\end{equation}

\hspace{0.6cm} In order to derive uniform estimates for ${v}_{\lambda \varepsilon}$, we recall first the following result.

\begin{lemma}[{\cite[Lemma 4.1]{CKK}}] Let $(u, v)$ be a solution of $(\ref{3.1})$, $(\ref{3.2})$. Then, the function
\begin{equation} t\mapsto {I}^{\lambda}_\varepsilon(u;v)(t)=\frac{1}{2\lambda}\left|[v(t)-{f}^*_\varepsilon(u(t))]^+\right|^2_H +\frac{1}{2\lambda}\left|[{f}_{*\varepsilon}(u(t))-v(t)]^+\right|^2_H \label{4.5}\end{equation}

is absolutely continuous on $[0,T]$ and
\begin{equation*} \frac{d}{dt} {I}_\varepsilon^{\lambda}(u;v)\leq  \langle \partial {I}_\varepsilon^{\lambda}(u;v),v'\rangle_H+L_0|u'|_H |\partial {I}_\varepsilon^{\lambda}(u;v)|_H   \label{}\end{equation*}

a.e. in $(0,T)$, where $L_0$ is a common Lipschitz constant of $f_*$ and $f^*$ $($see Lemma $2.1 (2))$. \qed
\end{lemma}

\hspace{0.6cm} Now, testing Eq. (\ref{3.2}) by $ v'_{\lambda \varepsilon}$ we see in view of Lemma 4.1 that
\begin{align} \nonumber\frac{1}{2}|v'_{\lambda \varepsilon}|_H^2&+\frac{d}{dt}\left(\frac{\kappa}{2}|\nabla v_{\lambda \varepsilon}|_H^2 +{I}^\lambda_\varepsilon(u_{\lambda\varepsilon};v_{\lambda\varepsilon})\right)\\ &\leq
L_0^2 |u'_{\lambda \varepsilon}|_H^2+\frac{1}{4}\left|\partial{I}^\lambda_\varepsilon(u_{\lambda\varepsilon};v_{\lambda\varepsilon})\right|_H^2+R_2
\label{4.6}\end{align}

 a.e. on $(0,T)$, where $R_2>0$ is a constant independent of $\varepsilon$. Then, testing  Eq. (\ref{3.2}) by $-\kappa\Delta v_{\lambda \varepsilon}$ yields
\begin{align} \nonumber \kappa^2|\Delta v_{\lambda \varepsilon}|_H^2+\frac{d}{dt}\left(\frac{\kappa}{2}|\nabla v_{\lambda \varepsilon}|_H^2\right) &\leq \left\langle \kappa\Delta v_{\lambda \varepsilon},\partial{I}^\lambda_\varepsilon(u_{\lambda\varepsilon};v_{\lambda\varepsilon})\right\rangle_H \\ &+\frac{1}{2}\kappa^2 |\Delta v_{\lambda \varepsilon}|_H^2+ 4R_2
\label{4.7}\end{align}

a.e. on $(0,T)$. We evaluate the first term on the right-hand side of this inequality as follows
\begin{align*} &\left\langle \partial{I}^\lambda_\varepsilon(u_{\lambda\varepsilon};v_{\lambda\varepsilon}), \kappa\Delta v_{\lambda \varepsilon}\right\rangle_H = \\
&\left\langle \frac{\kappa}{\lambda} \left[v_{\lambda \varepsilon}-{f}^*_\varepsilon(u_{\lambda\varepsilon})\right]^+,\Delta \left(v_{\lambda\varepsilon}-{f}^*_\varepsilon(u_{\lambda\varepsilon})\right)\right\rangle_H
 +\left\langle \frac{\kappa}{\lambda} \left[v_{\lambda \varepsilon}-{f}^*_\varepsilon(u_{\lambda\varepsilon})\right]^+,\Delta {f}^*_\varepsilon(u_{\lambda\varepsilon})\right\rangle_H +
\\
&\left\langle \frac{\kappa}{\lambda} \left[{f}_{*\varepsilon}(u_{\lambda\varepsilon})-v_{\lambda \varepsilon}\right]^+,\Delta \left({f}_{*\varepsilon}(u_{\lambda\varepsilon})-v_{\lambda\varepsilon}\right)\right\rangle_H
 +\left\langle \frac{\kappa}{\lambda} \left[{f}_{*\varepsilon}(u_{\lambda\varepsilon})-v_{\lambda \varepsilon}\right]^+,-\Delta {f}_{*\varepsilon}(u_{\lambda\varepsilon})\right\rangle_H
\\
&=-\frac{\kappa}{\lambda} \left|\nabla\left[v_{\lambda \varepsilon}-{f}^*_\varepsilon(u_{\lambda\varepsilon})\right]^+\right|_H^2
 +\left\langle \frac{\kappa}{\lambda} \left[v_{\lambda \varepsilon}-{f}^*_\varepsilon(u_{\lambda\varepsilon})\right]^+,\Delta {f}^*_\varepsilon(u_{\lambda\varepsilon})\right\rangle_H
\\
&\*\hspace{0.3cm}-\frac{\kappa}{\lambda} \left|\nabla\left[{f}_{*\varepsilon}(u_{\lambda\varepsilon})-v_{\lambda \varepsilon}\right]^+\right|_H^2
 +\left\langle \frac{\kappa}{\lambda} \left[{f}_{*\varepsilon}(u_{\lambda\varepsilon})-v_{\lambda \varepsilon}\right]^+,-\Delta {f}_{*\varepsilon}(u_{\lambda\varepsilon})\right\rangle_H
\\
&\leq\frac{1}{8\lambda^2} \left\{\left|\left[v_{\lambda \varepsilon}-{f}^*_\varepsilon(u_{\lambda\varepsilon})\right]^+\right|_H^2+\left|\left[{f}_{*\varepsilon}(u_{\lambda\varepsilon})-v_{\lambda \varepsilon}\right]^+\right|_H^2\right\} \\
&\*\hspace{4.6cm} +2\kappa^2\left(\left|\Delta {f}^*_\varepsilon(u_{\lambda\varepsilon})\right|_H^2 +\left|\Delta {f}_{*\varepsilon}(u_{\lambda\varepsilon})\right|_H^2\right)
\\
&\leq\frac{1}{8} \left|\partial{I}^\lambda_\varepsilon(u_{\lambda\varepsilon};v_{\lambda\varepsilon})\right|_H^2 +2\kappa^2\left(\left|\Delta {f}^*_\varepsilon(u_{\lambda\varepsilon})\right|_H^2 +\left|\Delta {f}_{*\varepsilon}(u_{\lambda\varepsilon})\right|_H^2\right).
  \label{}\end{align*}

Observing that $\Delta f(u)=f''(u)|\nabla u|_H^2+f'(u)\Delta u$ and invoking the Gagliardo-Nirenberg inequality from the last inequality we obtain
\begin{align*} &\left\langle \partial{I}^\lambda_\varepsilon(u_{\lambda\varepsilon};v_{\lambda\varepsilon}), \kappa\Delta v_{\lambda \varepsilon}\right\rangle_H
\leq\frac{1}{8} \left|\partial{I}^\lambda_\varepsilon(u_{\lambda\varepsilon};v_{\lambda\varepsilon})\right|_H^2 +2\kappa^2 R_3\left(|u_{\lambda\varepsilon}|_H^2+|\Delta u_{\lambda\varepsilon}|_H^2\right)
  \label{}\end{align*}

for a constant $R_3>0$ independent of $\varepsilon$, so that (\ref{4.7}) implies that
\begin{align} \nonumber \kappa^2|\Delta v_{\lambda \varepsilon}|_H^2+\frac{d}{dt}\left(\frac{\kappa}{2}|\nabla v_{\lambda \varepsilon}|_H^2\right) &\leq \frac{1}{8} \left|\partial{I}^\lambda_\varepsilon(u_{\lambda\varepsilon};v_{\lambda\varepsilon})\right|_H^2 \\ & + 2\kappa^2 R_3\left(|u_{\lambda\varepsilon}|_H^2+|\Delta u_{\lambda\varepsilon}|_H^2\right)
\label{4.8}\end{align}

a.e. on $(0,T)$. Similarly, testing Eq. (\ref{3.2}) by $\partial{I}^\lambda_\varepsilon(u_{\lambda\varepsilon};v_{\lambda\varepsilon})$ and using Lemma 3.1 we see that
\begin{align} \nonumber \frac{1}{2} \left|\partial{I}^\lambda_\varepsilon(u_{\lambda\varepsilon};v_{\lambda\varepsilon})\right|_H^2
&+ \frac{d}{dt} {I}^\lambda_\varepsilon(u_{\lambda\varepsilon};v_{\lambda\varepsilon})\leq L_0^2|u'_{\lambda \varepsilon}|_H^2\\
&+2\kappa^2 R_3\left(|u_{\lambda\varepsilon}|_H^2+|\Delta u_{\lambda\varepsilon}|_H^2\right)+4R_2
\label{4.9}\end{align}

a.e. on $(0,T)$.

\hspace{0.6cm}  Next, taking the sum of the inequalities $(\ref{4.6})$, $(\ref{4.8})$, and $(\ref{4.9})$, then integrating the result from $0$ to $T$ we see in view of  $(\ref{3.4})$, $(\ref{4.1})$, and $(\ref{4.4})$ that
\begin{align} 
|v_{\lambda \varepsilon}'|_{L^2(0,T;H)}+\kappa|\Delta v_{\lambda \varepsilon}|_{L^2(0,T;H)}
+\left|\partial{I}^\lambda_\varepsilon(u_{\lambda\varepsilon};v_{\lambda\varepsilon})\right|_{L^2(0,T;H)}\leq R_4 \label{4.10}
\end{align}

for a constant $R_4>0$ independent of $\varepsilon$.

\hspace{0.6cm} In view of $(\ref{4.10})$, similarly to $(\ref{3.16})$ we obtain
\begin{equation*} \frac{\kappa}{2}\int_0^T|\nabla v_{\lambda \varepsilon}(\tau)|^2_H d\tau\leq R_5,  \label{}\end{equation*}

where $R_5>0$ is a constant independent of $\varepsilon$.  Invoking the Poincar\'{e} inequality  we further derive
\begin{equation} \frac{\kappa}{2}\int_0^T| v_{\lambda \varepsilon}(\tau)|^2_H d\tau\leq R_5 C_P. \label{4.13}\end{equation}

Then, from $(\ref{4.5})$, $(\ref{4.4})$, and $(\ref{4.13})$ we see that
\begin{equation} \int_0^T I_\lambda^{\varepsilon}(u_{\lambda \varepsilon}(\tau),v_{\lambda \varepsilon}(\tau)) d\tau \leq \frac{1}{\lambda}R_6, \label{4.14}\end{equation}

where $R_6>0$ is a constant independent of $\varepsilon$ and $\lambda$. Taking account of $(\ref{4.6})$, $(\ref{4.8})$, $(\ref{4.9})$, $(\ref{4.1})$, $(\ref{4.5})$,  $(\ref{4.14})$ and reasoning as in derivation of $(\ref{3.17})$ we infer that
\begin{equation} |\nabla v_{\lambda \varepsilon_n}|_{L^\infty(0,T:H)}\leq R(\lambda) \label{4.15}\end{equation}

for a positive constant  $R(\lambda)$ dependent on $\lambda$, $\kappa$, $L_0$, $T$, $R_1$, $R_2$, $R_3$, but independent of $\varepsilon$.

\hspace{0.6cm} As before, the estimates $(\ref{4.10})$, $(\ref{4.15})$ imply the existence of  a null sequence  $\{\varepsilon_n\}_{n\geq 1}$  and a function $v_\lambda$ such that
\begin{align}
& \left. \begin{array}{lc} & {v}_{\lambda \varepsilon_n} \to v_\lambda \quad
\mbox{weakly in } W^{1,2}(0,T; H) \cap
L^2(0,T;H^2(\Omega))  \vspace{0.1cm}\\
& \*\hspace{3cm} \mbox{and weakly-star in } L^{\infty}(0,T; V).
\end{array}
\right. \label{4.16}
\end{align}

In particular, we also have
\begin{equation} {v}_{\lambda \varepsilon_n} \to v_\lambda \quad  \mbox{in } C([0,T];H).\label{4.17}\end{equation}

\hspace{0.6cm} Now, from the convergences  $(\ref{4.2})$, $(\ref{4.3})$, $(\ref{4.16})$, $(\ref{4.17})$, Lemma 2.1 (2), (3), and the Arzela--Ascoli theorem we see that the pair $\{u_\lambda,v_\lambda\}$, $\lambda>0$, is a solution of the following system, which we denote by $(P)_\lambda$:
\begin{equation}
 u' -\Delta u =h(u,v) \quad \text{in } H  \; \text{a.e. on } [0,T]
,\label{4.18}
\end{equation}
\begin{equation} v'-\kappa\Delta v+ \partial {I}^{\lambda}(u;v)= g(u,v)   \quad \text{in } H  \; \text{a.e. on } [0,T] , \label{4.19}\end{equation}
\begin{equation}u =v= 0 \quad \text{on } \partial\Omega \;\; \text{a.e. on } [0,T], \label{4.20}\end{equation}
\begin{equation}u(0) =
u(T), \quad v(0) = v(T) \quad \text{in } H ,
\label{4.21}\end{equation}

where $\partial {I}^{\lambda}(u;v)$ is defined  in $(\ref{2.2})$.

\hspace{0.6cm} A solution to
$(P)_{\lambda}$ is a pair of functions $\{u,v\}$ such that $u, v\in W^{1,2}(0,T;H)\cap L^\infty(0,T;V)\cap
L^2(0,T;H^2 (\Omega ))$ and  $(\ref{4.18})$--$(\ref{4.21})$ hold.

\hspace{0.6cm}  We note that the validity of the periodic condition $(\ref{4.21})$ follows from $(\ref{3.4})$ and $(\ref{4.3})$, $(\ref{4.17})$.

\subsection{Passage-to-the-limit: $\lambda\to 0$}

 We now derive a priori estimates uniform with respect to the parameter $\lambda>0$ for solutions $(u_{\lambda}, v_{\lambda})$ of Problem $(P)_{\lambda}$.

\hspace{0.6cm}  To this aim, we note that the constants $R_1$, $R_4$, and $R_5$ in the uniform estimates of the previous subsection do not depend on $\lambda$. Hence, for a solution $(u_\lambda,v_\lambda)$ of Problem $(P_\lambda)$ from $(\ref{4.1})$, $(\ref{4.10})$, and $(\ref{4.13})$ we have
\begin{align*} 
|u_{\lambda}'|_{L^2(0,T;H)}&+|\Delta u_{\lambda}|_{L^2(0,T;H)}
+|\nabla u_{\lambda}|_{L^\infty(0,T;H)}\\ &+|v_{\lambda }'|_{L^2(0,T;H)}+\kappa|\Delta v_{\lambda }|_{L^2(0,T;H)}
+\left|\partial{I}^\lambda(u_{\lambda};v_{\lambda})\right|_{L^2(0,T;H)}\leq R_7
\end{align*}

and
\begin{equation}|v_{\lambda }|_{L^2(0,T;H)}\leq R_7, \label{4.22}\end{equation}

for a constant $R_7>0$ independent of $\lambda$. In particular, as above we conclude that there exists a null sequence
$\lambda_n$, $n\geq 1$, in $(0,1]$ and functions $u, v$ such that
\begin{align}
& \left. \begin{array}{lc} &\*\hspace{-1cm} u_{\lambda_n} \to u  \quad
  \vspace{0.1cm}\mbox{weakly in }  W^{1,2}(0,T; H) \cap
L^2(0,T;H^2(\Omega))\\
& \*\hspace{1.8cm} \mbox{and weakly-star in } L^{\infty}(0,T; V)
  \vspace{0.1cm}\\
& \*\hspace{2.3cm}\mbox{and, thus, strongly in } C([0,T];H),
\end{array}
\right. \label{4.23}\\
&\*\hspace{-0.1cm}v_{\lambda_n} \to v\quad \mbox{weakly in }
W^{1,2}(0,T; H) ,\label{4.24}
\\
&\*\hspace{-0.6cm}\kappa\Delta v_{\lambda_n} \to \kappa\Delta v\quad \mbox{weakly in }
L^{2}(0,T; H) ,\label{4.25}
\\
&\*\hspace{-0.6cm}\partial {I}^{\lambda_n}(u_n,v_n) \to \xi\quad \mbox{weakly in }
L^{2}(0,T; H) \label{4.26}
\end{align}

for some function $\xi\in L^{2}(0,T; H)$.

\hspace{0.6cm}  Below, we show that along with the convergences $(\ref{4.23})$--$(\ref{4.26})$ we also have
\begin{equation} v_{\lambda_n} \to
v \quad \mbox{strongly in } C([0,T];H).\label{4.27}\end{equation}

To this end, take two arbitrary $i,j\geq 1$ with $i\neq j$ and denote $u_i:=u_{\lambda_i}$, $v_i:=v_{\lambda_i}$. Then,  from $(\ref{4.19})$ it follows that
$$v'_{j}-v'_{i}-\kappa(\Delta v_{j}-\Delta v_{i})+\partial {I}^{\lambda_j}(u_j;v_j)-\partial {I}^{\lambda_i}(u_i;v_i)=g(u_j,v_j)-g(u_i,v_i),$$

Testing this equality by $v_j-v_i$, using the Lipschitz continuity
of $g$ and invoking Young's inequality we have
\begin{align}\nonumber \frac{1}{2}\frac{d}{dt}|v_{j}-v_{i}|_H^2&+\kappa|\nabla(v_{j}-v_{i})|_H^2+\langle\partial {I}^{\lambda_j}(u_j;v_j)-\partial {I}^{\lambda_i}(u_i;v_i),v_j-v_i\rangle_H \\ & \leq L_g|u_j-u_i|_H|v_j-v_i|_H+L_*|v_j-v_i|^2_H,\label{4.28}\end{align}

where $L_g$ and $L_*$ are the Lipschitz constants of $g$ in $u$ and $v$, respectively.

Setting
\begin{align}S_{ij}^\lambda=\langle\partial {I}^{\lambda_j}(u_j;v_j)-\partial {I}^{\lambda_i}(u_i;v_i),v_j-v_i \rangle_H.\label{4.29}\end{align}

we see from $(\ref{2.2})$  that
\begin{align*}
S_{ij}^\lambda=\bigg\langle&\frac{1}{\lambda_j}[v_j-{f}^*(u_j)]^+-\frac{1}{\lambda_j}[{f}_*(u_j)-v_j]^+\\-&\frac{1}{\lambda_i}[v_i-{f}^*(u_i)]^++\frac{1}{\lambda_i}[{f}_*(u_i)-v_i]^+,v_j-v_i\bigg\rangle_H.
\end{align*}

We have nine
possible cases to estimate the value of $S_{ij}^\lambda$ from below.  First, assuming that $v_j\geq
{f}^*(u_j)$, $v_i\geq {f}^*(u_i)$ we obtain
\begin{align*}
S_{ij}^\lambda=\bigg\langle\langle&\frac{1}{\lambda_j}(v_j-{f}^*(u_j))-\frac{1}{\lambda_i}(v_i-{f}^*(u_i)),\\
&\lambda_j\frac{1}{\lambda_j}(v_j-{f}^*(u_j))-\lambda_i\frac{1}{\lambda_i}(v_i-{f}^*(u_i))+{f}^*(u_j)-{f}^*(u_i)\bigg\rangle_H
\\ \geq \lambda_j& |\partial {I}^{\lambda_j}(u_j;v_j)|_H^2+\lambda_i |\partial
{I}^{\lambda_i}(u_i;v_i)|_H^2-(\lambda_j+\lambda_i) |\partial
{I}^{\lambda_j}(u_j;v_j)|_H |\partial {I}^{\lambda_i}(u_i;v_i)|_H\\
&- (|\partial {I}^{\lambda_j}(u_j;v_j)|_H +|\partial
{I}^{\lambda_i}(u_i;v_i)|_H)|{f}^*(u_j)-{f}^*(u_i)|_H.\end{align*}

Second, when $v_j\geq {f}^*(u_j)$, ${f}_*(u_i)\leq v_i< {f}^*(u_i)$ we
see that
\begin{align*}
S_{ij}^\lambda=\bigg\langle\frac{1}{\lambda_j}(v_j-{f}^*(u_j)),v_j-v_i\bigg\rangle_H \geq
-|\partial {I}^{\lambda_j}(u_j;v_j)|_H
|{f}^*(u_j)-{f}^*(u_i)|_H.\end{align*}

Third, if $v_j\geq {f}^*(u_j)$, $ v_i< {f}_*(u_i)$, then
\begin{align*}
S_{ij}^\lambda&=\bigg\langle\frac{1}{\lambda_j}(v_j-{f}^*(u_j))+\frac{1}{\lambda_i}({f}_*(u_i)-v_i),v_j-v_i\bigg\rangle_H
\\ &\geq -\left(|\partial {I}^{\lambda_j}(u_j;v_j)|_H+|\partial
{I}^{\lambda_i}(u_i;v_i)|_H\right) |{f}^*(u_j)-{f}^*(u_i)|_H.\end{align*}

The reasoning in the remaining cases:
\begin{align*}
{f}_*(u_j)< v_j&< {f}^*(u_j), \; v_i\geq {f}^*(u_i)\\
{f}_*(u_j)< v_j&< {f}^*(u_j), \; {f}_*(u_i)\leq v_i< {f}^*(u_i)\\
{f}_*(u_j)< v_j&< {f}^*(u_j), \; v_i< {f}_*(u_i) \\
v_j&\leq {f}_*(u_j), \; v_i\geq {f}^*(u_i)\\
v_j&\leq {f}_*(u_j), \; {f}_*(u_i)\leq v_i< {f}^*(u_i)\\
v_j&\leq {f}_*(u_j), \; v_i< {f}_*(u_i) \\
\end{align*}

is fully symmetric and is left to the reader. Consequently, we always have
\begin{align} \nonumber
S_{ij}^\lambda&\geq -(\lambda_j+\lambda_i) |\partial
{I}^{\lambda_j}(u_j;v_j)|_H |\partial {I}^{\lambda_i}(u_i;v_i)|_H\\  \nonumber
&- (|\partial {I}^{\lambda_j}(u_j;v_j)|_H +|\partial
{I}^{\lambda_i}(u_i;v_i)|_H)\left(|{f}^*(u_j)-{f}^*(u_i)|_H+|{f}_*(u_j)-{f}_*(u_i)|_H\right)\\
& =: -\delta_{ij}^\lambda. \label{4.30}\end{align}

\hspace{0.6cm} For convenience, denote $\bar{u}:=u_j-u_i$, $\bar{v}:=v_j-v_i$. Then, from $(\ref{4.28})$--$(\ref{4.30})$ we infer that
\begin{align*}\frac{1}{2}\frac{d}{dt}|\bar{v}|_H^2+\kappa|\nabla \bar{v}|_H^2\leq L_g|\bar{u}|_H|\bar{v}|_H+L_*|\bar{v}|^2_H+\delta_{ij}^\lambda.\label{}\end{align*}

Invoking the Poincar\'{e} inequality, from this inequality we obtain
\begin{align*}\frac{1}{2}\frac{d}{dt}|\bar{v}|_H^2\leq -\left(\frac{\kappa}{C_P}-L_*\right)|\bar{v}|^2_H+ L_g|\bar{u}|_H|\bar{v}|_H+\delta_{ij}^\lambda.\label{}\end{align*}

Denote $\displaystyle c_0:=\frac{\kappa}{C_P}-L_*$ ($c_0>0$ by Hypothesis $(H1)$). We further have
\begin{align*}\frac{1}{2}\frac{d}{dt}\left(e^{2c_0 t}|\bar{v}|_H^2\right)\leq 2 e^{2c_0 t}\left(L_g|\bar{u}|_H|\bar{v}|_H+\delta_{ij}^\lambda\right).\label{}\end{align*}

Integrating this inequality from $0$ to $t\in (0,T]$ and using H\"{o}lder's inequality we see in view of  $(\ref{4.22})$ that

\begin{align}\nonumber |\bar{v}(t)|_H^2 &\leq e^{-2c_0 t}|\bar{v}(0)|_H^2+2\int_0^t \left(L_g|\bar{u}(\tau)|_H|\bar{v}(\tau)|_H+\delta_{ij}^\lambda(\tau)\right) d\tau \\
&\leq e^{-2c_0 t}|\bar{v}(0)|_H^2+4L_g R_7|\bar{u}|_{L^2(0,T;H)}+2\int_0^t \delta_{ij}^\lambda(\tau) d\tau.\label{4.31}\end{align}

Taking $t=T$ in $(\ref{4.31})$ and substituting the resulting inequality, using the fact that $\bar{v}(T)=\bar{v}(0)$, into $(\ref{4.31})$ we obtain
\begin{align*} |\bar{v}(t)|_H^2 \leq \frac{2e^{2c_0 T}}{e^{2c_0 T}-1}\left(2L_g R_7|\bar{u}|_{L^2(0,T;H)}+\int_0^T \delta_{ij}^\lambda(\tau) d\tau\right).\label{}\end{align*}

Applying Gronwall's inequality to this inequality  we conclude in view of  the convergence $(\ref{4.23})$ that $v_i$, $i\geq 1$, is a Cauchy sequence in the space $C([0,T];H)$. Hence, according
to $(\ref{4.24})$ we obtain the convergence $(\ref{4.27})$.

\hspace{0.6cm} To show  that $v\in L^\infty(0,T;V)$, we note that $(\ref{4.19})$ implies that
\begin{equation*}
g(u(t),v(t))- v'(t)-\partial I^\lambda(u(t),v(t))+ v(t)\in \partial\psi(v(t)) \quad\mbox{for a.e.} \; t\in [0,T], \label{}
\end{equation*}

where
\begin{equation*}\psi(v):= \left\{
\begin{array}{ll}
\frac{1}{2}|v|_H^2 +\frac{1}{2}|\nabla v|_H^2 & \mbox{if} \;v \in V ,\\
+\infty  & \mbox{if} \;v \in H\setminus V,
\end{array}\right. \quad \partial\varphi(v)=v-\Delta v.
\label{}\end{equation*}

Since all the functions on the left-hand side of the inclusion above
belong to $L^2(0,T;H)$, from \cite[Lemma
3.3]{Brezis} it follows that the function
$$t\to \varphi(v(t))=\frac{1}{2}|v(t)|_H^2
+\frac{1}{2}|\nabla v(t)|_H^2$$ is absolutely continuous. Consequently, the function
$
t\to |v(t)|_V=(|v(t)|_H^2
+|\nabla v(t)|_H^2)^{1/2}
$
is continuous. Hence, the set $\{v(t): \; t\in [0,T]\}$ is
bounded in the space $V$. Since the embedding $H
\hookrightarrow V'$ is dense, the function  $t\to v(t)$ which is
continuous from $[0,T]$ to $H$ is also continuous from $[0,T]$ to
$\omega$-$V$. Since in a reflexive Banach space the
weak convergence  coupled with the convergence of norms implies
the norm convergence, we conclude that the  function
$t\to v(t)$ is continuous from $[0,T]$ to $V$. In particular, $v\in L^\infty(0,T;V)$.

\hspace{0.6cm} Given the convergences
$(\ref{4.23})$--$(\ref{4.27})$, to finish the proof that the pair $\{u,v\}$ is a
solution to Problem  $(P)$ it remains to show that
\begin{equation} \xi\in \partial I(u;v) \quad    \hbox{a.e. on} \
(0,T).\label{4.32}\end{equation}

To this end, let $z$ be an arbitrary
function from $L^2(0,T;H)$ such that $z\in[{f}_*(u),{f}^*(u)]$  a.e.
on $Q$. For every $n\geq 1$, define $z_n$ to be the pointwise projection of $z$ onto the set $[{f}_*(u_n),{f}^*(u_n)].$
Then, $z_n\in [{f}_*(u_n),{f}^*(u_n)]$ a.e. on $Q$, $n\geq 1$, and
$z_n\to z$ in $L^2(0,T;H)$ as $n\to\infty$. Consequently, since the operator $\partial {I}^{\lambda_n}(u_n;\cdot)$ is the subdifferential of the function ${I}^{\lambda_n}(u_n;\cdot)$, from the definition of subdifferential we have
\begin{equation} \langle\partial {I}^{\lambda_n}(u_n;v_n), z_n-v_n\rangle_H \leq  {I}^{\lambda_n}(u_n;z_n)-{I}^{\lambda_n}(u_n;v_n)\leq -{I}^{\lambda_n}(u_n;v_n)\leq 0,   \label{4.33}\end{equation}

$ n\geq 1.$ On the other hand, from the uniform boundedness of $\{\partial I^{\lambda_n}(u_n;v_n)\}$, $n\geq 1$, in $L^2(0,T;H)$ and $(\ref{2.2})$ we see that
$$[v_n-{f}^*(u_n)]^+ +[{f}_*(u_n)-v_n]^+=\lambda_n \left|\partial {I}^{\lambda_n}(u_n;v_n)\right|\to 0$$

in $L^2(0,T;H)$ as $n\to\infty$. Therefore, we infer that $v\in[{f}_*(u), {f}^*(u)]$  a.e. on $Q$. Passing now to the limit as
$n\to\infty$  in $(\ref{4.33})$ we conclude that $(\ref{4.32})$ holds
and $\{u,w\}$ is thus a solution to problem $(P)$.

\hspace{0.6cm} Finally, we note that the periodicity condition $(iv)$ in Definition 2.1 follows from $(\ref{4.21})$, $(\ref{4.23})$, and $(\ref{4.27})$.

\vspace{0.5cm}

\textbf{Acknowledgements.}  The authors want to thank the
anonymous referees for their valuable suggestions and remarks
which helped to improve the manuscript.

\bibliographystyle{amsplain}

\begin{thebibliography}{40}





\bibitem {BS} M. Brokate and J. Sprekels, \emph{Hysteresis and
    Phase
    Transitions}.  Appl. Math. Sci., 121, Springer-Verlag,
    New York, 1996.


\bibitem {Krejci} P. Krej\v{c}\'{\i}, \emph{Hysteresis,
    Convexity and
    Dissipation in Hyperbolic Equations}. Gakuto Int. Ser. Math.
    Sci. Appl., Vol.~8, Gakk\={o}tosho, Tokyo, 1996.



\bibitem{Visintin} A. Visintin,
\emph{Differential Models of Hysteresis}. Appl. Math. Sci. 111,
Springer-Verlag, Berlin, 1994.






\bibitem{KKV} N. Kenmochi, T. Koyama and A. Visintin,
\emph{On a class of variational inequalities with memory terms}, in:
Progress in partial differential equations: elliptic and parabolic
problems. Pitman Res. Notes Math., Ser. 266, Longman, Harlow (1992),
164--175.



\bibitem{KeV} N. Kenmochi and A. Visintin,
\emph{Asymptotic stability for parabolic variational inequalities
with hysteresis}, in: Models of hysteresis. Pitman Res. Notes Math.,
Ser. 286, Longman, Harlow (1993), 59--70.



\bibitem{KeKM} N. Kenmochi, T. Koyama and G.H. Meyer,
    \emph{Parabolic PDEs with hysteresis and quasivariational
    inequalities}. Nonlinear Anal. 34 (1998),
    665--686.


\bibitem{CKK} P. Colli, N. Kenmochi and M. Kubo, \emph{A phase
    field model with temperature dependent constraint}. J.
    Math. Anal. Appl. 256 (2001), 668--685.




\bibitem{KMO} N. Kenmochi, E. Minchev and T. Okazaki, \emph{On
    a system of nonlinear PDE's with diffusion and  hysteresis
    effects}. Adv. Math. Sci. Appl.  14 (2004),
    633--664.



\bibitem{Kubo} M. Kubo, \emph{A filtration model with hysteresis}. J. Differential Equations  2001 (2004), no. 1, 75--98.


\bibitem{KTT} P. Krej\v{c}\'{\i}, A.A. Tolstonogov and S.A.
    Timoshin, \emph{A control problem in phase transition
    modeling}. NoDEA, Nonlinear Differ. Equ. Appl. 22, No. 4
    (2015), 513--542.


\bibitem{AT} T. Aiki and S.A. Timoshin, \emph{Existence and uniqueness for a concrete carbonation process with hysteresis}, J.
    Math. Anal. Appl. 449 (2017), no. 2, 1502-1519.




\bibitem{Krejci1} P. Krej\v{c}\'{\i}, \emph{Hysteresis and periodic solutions of semilinear and quasilinear wave equations}, Math. Z. 193 (1986), no. 2, 247–264.


\bibitem{Krejci2} P. Krej\v{c}\'{\i}, \emph{Forced periodic vibrations of an elastic system with elastico-plastic damping}, 
Apl. Mat. 33, No. 2, 145-153 (1988).


\bibitem{Krejci3} P. Krej\v{c}\'{\i},  \emph{Reliable solutions to the problem of periodic oscillations of an elastoplastic beam}, Internat. J. Non-Linear Mech. 37 (2002), no. 8, 1337–1349.



\bibitem{GHM} I.G. G\"{o}tz, K.-H. Hoffmann, and A.M. Meirmanov, \emph{Periodic solutions of the Stefan problem with hysteresis-type boundary conditions}, Manuscripta Math. 78 (1993), no. 2, 179–-199.




\bibitem{comwa} J. Zheng, Y. Ke, and Y. Wang, \emph{Periodic solutions to a heat equation with hysteresis in the source term}, Comput. Math. Appl. 69 (2015), no. 2, 134–143.



\bibitem{Quittner} P. Quittner, \emph{Liouville theorems, universal estimates and periodic solutions for cooperative parabolic Lotka-Volterra systems}.  J. Differential Equations 260 (2016), no. 4, 3524--3537.


\bibitem{Yamada} Y. Yamada, \emph{Periodic solutions of certain nonlinear parabolic differential equations in domains with periodically moving boundaries}. Nagoya Math. J. 70 (1978), 111–-123.


\bibitem{Kenmochi}  N. Kenmochi,  \emph{Solvability of nonlinear evolution equations with
time-dependent constraints and applications}, Bull. Fac. Educ.,
Chiba Univ., Part 2, 30 (1981), 1--87.


\bibitem{Brezis} H. Brezis,
    \emph{Op\'{e}rateurs maximaux monotones et semi-groupes de contractions dans les espaces de Hilbert}. North-Holland, Amsterdam, London, 1973





\end{thebibliography}

\end{document}